\documentclass[12pt]{amsart}
\usepackage{graphicx}
\usepackage{graphics}
\usepackage{amsmath}
\usepackage{amscd}
\usepackage{latexsym}
\usepackage{subfigure}
\usepackage{caption}

\begin{document}

\textwidth 5.9in
\textheight 7.9in

\evensidemargin .75in
\oddsidemargin .75in

\newtheorem{Thm}{Theorem}
\newtheorem{Lem}[Thm]{Lemma}
\newtheorem{Cor}[Thm]{Corollary}
\newtheorem{Prop}[Thm]{Proposition}
\newtheorem{Rm}{Remark}

\def\a{{\mathbb a}}
\def\C{{\mathbb C}}
\def\A{{\mathbb A}}
\def\B{{\mathbb B}}
\def\D{{\mathbb D}}
\def\E{{\mathbb E}}
\def\R{{\mathbb R}}
\def\P{{\mathbb P}}
\def\S{{\mathbb S}}
\def\Z{{\mathbb Z}}
\def\O{{\mathbb O}}
\def\H{{\mathbb H}}
\def\V{{\mathbb V}}
\def\Q{{\mathbb Q}}
\def\Cn{${\mathcal C}_n$}
\def\CM{\mathcal M}
\def\CG{\mathcal G}
\def\CH{\mathcal H}
\def\CT{\mathcal T}
\def\CF{\mathcal F}
\def\CA{\mathcal A}
\def\CB{\mathcal B}
\def\CD{\mathcal D}
\def\CP{\mathcal P}
\def\CS{\mathcal S}
\def\CZ{\mathcal Z}
\def\CE{\mathcal E}
\def\CL{\mathcal L}
\def\CV{\mathcal V}
\def\CW{\mathcal W}
\def\IC{\mathbb C}
\def\IF{\mathbb F}
\def\IK{\mathcal K}
\def\IL{\mathcal L}
\def\IP{\bf P}
\def\IR{\mathbb R}
\def\IZ{\mathbb Z}

\title{Topology of multiple log transforms of $4$-manifolds }
\author{Selman Akbulut}
\thanks{Partially supported by NSF grants DMS 0905917, DMS 1065955, and MPIM}
\address{Department  of Mathematics, Michigan State University,  MI, 48824}
\email{akbulut@math.msu.edu }
\date{\today}
\begin{abstract} 
Given a  $4$-manifold $X$ and an imbedding $T^{2}\times B^2 \subset X$, we describe an algorithm $X \mapsto X_{p,q}$ for drawing the handlebody of the $4$-manifold obtained from $X$ by $(p,q)$-logarithmic transforms along the parallel tori. By using this algorithm, we  obtain a simple handle picture of the Dolgachev surface $E(1)_{p,q}$, from that deduce that the exotic copy $E(1)_{p,q}\# \; 5 \bar{\C\P^2}$ of $E(1) \#\; 5 \bar{\C\P^2}$ differs from the original one by a  codimension zero simply connected Stein submanifold $M_{p,q}$, which are therefore examples of infinitely many Stein manifolds that are exotic copies of each other (rel boundaries). Furthermore, by a similar method we produce infinitely many simply connected Stein manifolds  $Z_{p}\subset E(1)_{p,2}\# \; 2 \bar{\C\P^2}$ with the same boundary and $b_{2}(Z_{p})=2$,  which are (absolutely) exotic copies of each other; this provides an alternative proof of a recent theorem of the author and Yasui \cite{ay4}.
Also, by using the description of $S^2 \times S^2$ as a union of two cusps glued along their boundaries, and by using this algorithm, we show that multiple log transforms along the tori in these cusps do not change smooth structure of $S^{2}\times S^{2}$.
\end{abstract}

\date{}
\maketitle

\setcounter{section}{-1}

\vspace{-.25in}

\section{Introduction}

Given a smooth $4$-manifold $X^4$ and an imbedding $T^2\times B^2\subset X^4$, the operation of removing this  $T^2\times B^2$ and then gluing it back by a nontrivial diffeomorphism of its boundary $\varphi: T^3\to T^3$  is called $T^2$-surgery operation. It is known that any smooth $4$-manifold can be obtained from $S^2\times S^2$ or $\C\P^2$ be a sequence of these operations
 \cite{i}. 
 
 \vspace{.05in}
 
A special version of this operation is called  {\it p-log transform}, where $ \varphi_{p}: T^{2}\times \partial B^2 \to T^{2}\times \partial B^2$  is the self-diffeomorphism  given by ($p\in \Z) $
\begin{equation*}
\varphi_{p}= \left(
\begin{array}{ccc}
1 &0 &0  \\
0 &p &1  \\
0 &1 &0
\end{array}
\right)
\end{equation*}

\begin{figure}[ht]  \begin{center}  
\includegraphics[width=.62\textwidth]{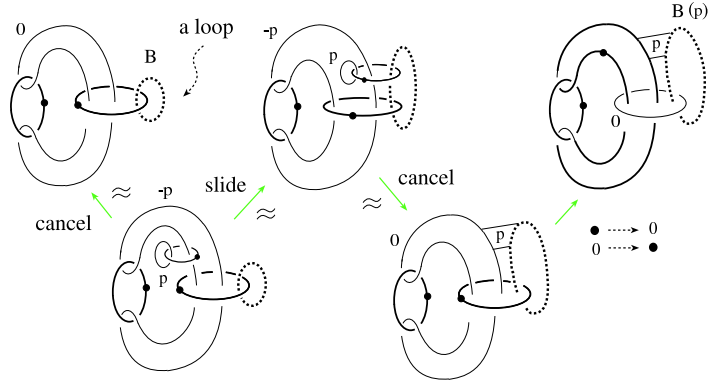}   
\caption{$p$ log-transform operation} 
\label{fig1}
\end{center}
\end{figure}

\noindent This operation is sometimes called  {\it  logarithmic transform of order p}.  Figure~\ref{fig1}describes this as a handlebody operation (c.f. \cite{a1} \cite{ay1}, \cite{g}).  This figure gives a recipe of how to modify a $4$-manifold handlebody, containing a framed torus $T^2\times D^2$,  to get the handlebody of the $p$-log transformed $4$-manifold along this torus. For example Figure~\ref{fig1} describes how the linking loop B is changed by this operation. 

The affect of this operation on complex surfaces had been studied by Kodaira and Dolgachev. The special class of complex surfaces  $E(1)_{p,q}$, which are obtained by performing  two log transforms of order $p, q \geq 2$ (relatively prime) on two parallel torus fibers of the elliptic surface $E(1)= \C\P^2\# 9\bar{\C\P}^2$ are called {\it Dolgachev Surfaces}. Donaldson proved, as smooth manifolds, they are exotic copies of $E(1)$ \cite{d}. 

\vspace{.05in}

The  algorithm above provides a useful tool understanding the change of the handlebody structure of the underlying manifold under p-log transform operation, especially when done multiple times.  For example as an exercise, the reader can check from the pictures that doing p-log transform operation to a torus followed  by q-log transform operation on the resulting torus is not the pq-log transform operation of the initial torus. One should be extra careful especially when applying multiple log transform operations one after another, since ordering of the circle factors of $T^2\times \partial B^2$ become important (note that different choices alter the matrix of $\varphi_{p}$).  Here we adapt one standard convention for (p,q)-log transform operation on parallel tori, where this choice is done from very beginning by fixing handlebody pictures of the parallel tori (Figure~\ref{fig5}). Previously in \cite{g} a  handle picture of $E(1)_{p,q}$ was given, which is unfortunately crowded by the $1$-handles drawn in the ``pair of balls'' notation. The approach taken in this paper (i.e. starting with parallel tori first,  then building the rest of the manifold) offers an easier alternative, enabling us to derive the desired corollaries discussed in this paper. Here I would like to thank K. Yasui for constructive comments.

\section{The algorithm $X\mapsto X_{p,q}$}
\begin{figure}[ht]  \begin{center}  
\includegraphics[width=.9\textwidth]{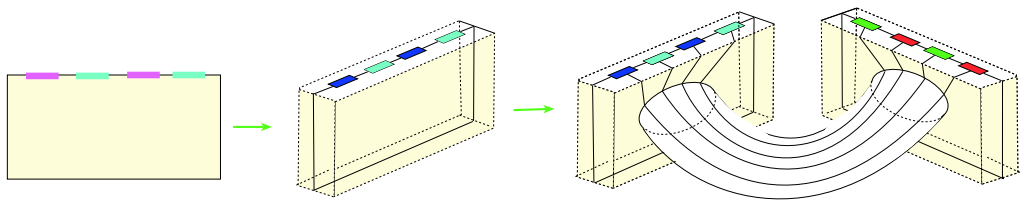}   
\caption{First  thicken $T^2$ to $T^2\times [0,1]$, then decompose it as  $T^2\times [0,1/2] \cup T^{2}\times [1/2,1]$ by face identification} 
\label{fig2}
\end{center}
\end{figure}

\begin{figure}[ht]  \begin{center}  
\includegraphics[width=.5\textwidth]{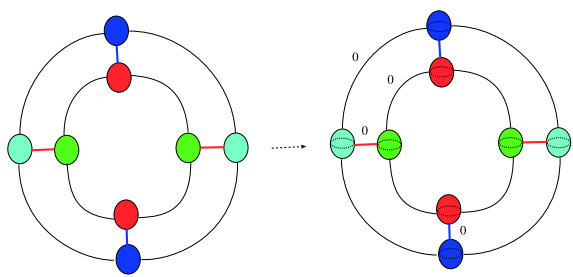}   
\caption{Thickening $T^2\times [0,1]$ to $T^{2}\times B^2$} 
 \label{fig3}
\end{center}
\end{figure}

\begin{figure}[ht]  \begin{center}  
\includegraphics[width=.8\textwidth]{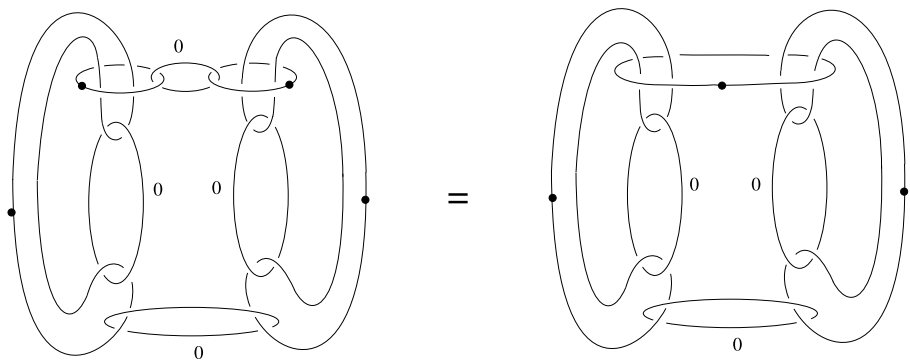}   
\caption{A handlebody picture of  $T^2\times B^2$ described as two parallel copies of $T^2 \times B^2$ glued side by side} 
 \label{fig4}
\end{center}
\end{figure}

\begin{figure}[ht]  \begin{center}  
\includegraphics[width=.5\textwidth]{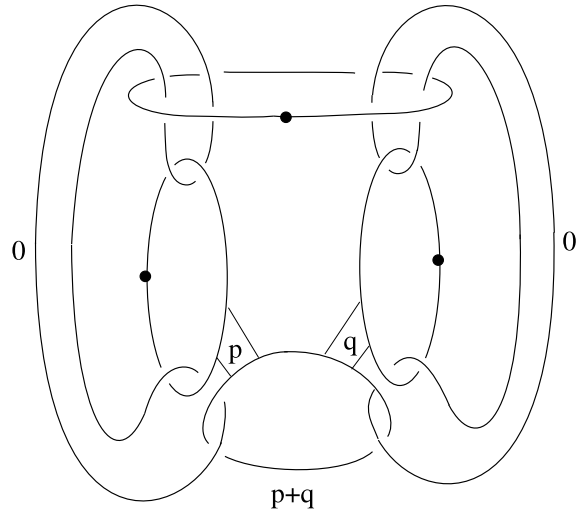}   
\caption{$(T^2\times B^2)_{p,q} $} 
 \label{fig5}
\end{center}
\end{figure}

\begin{figure}[ht]  \begin{center}  
\includegraphics[width=.8\textwidth]{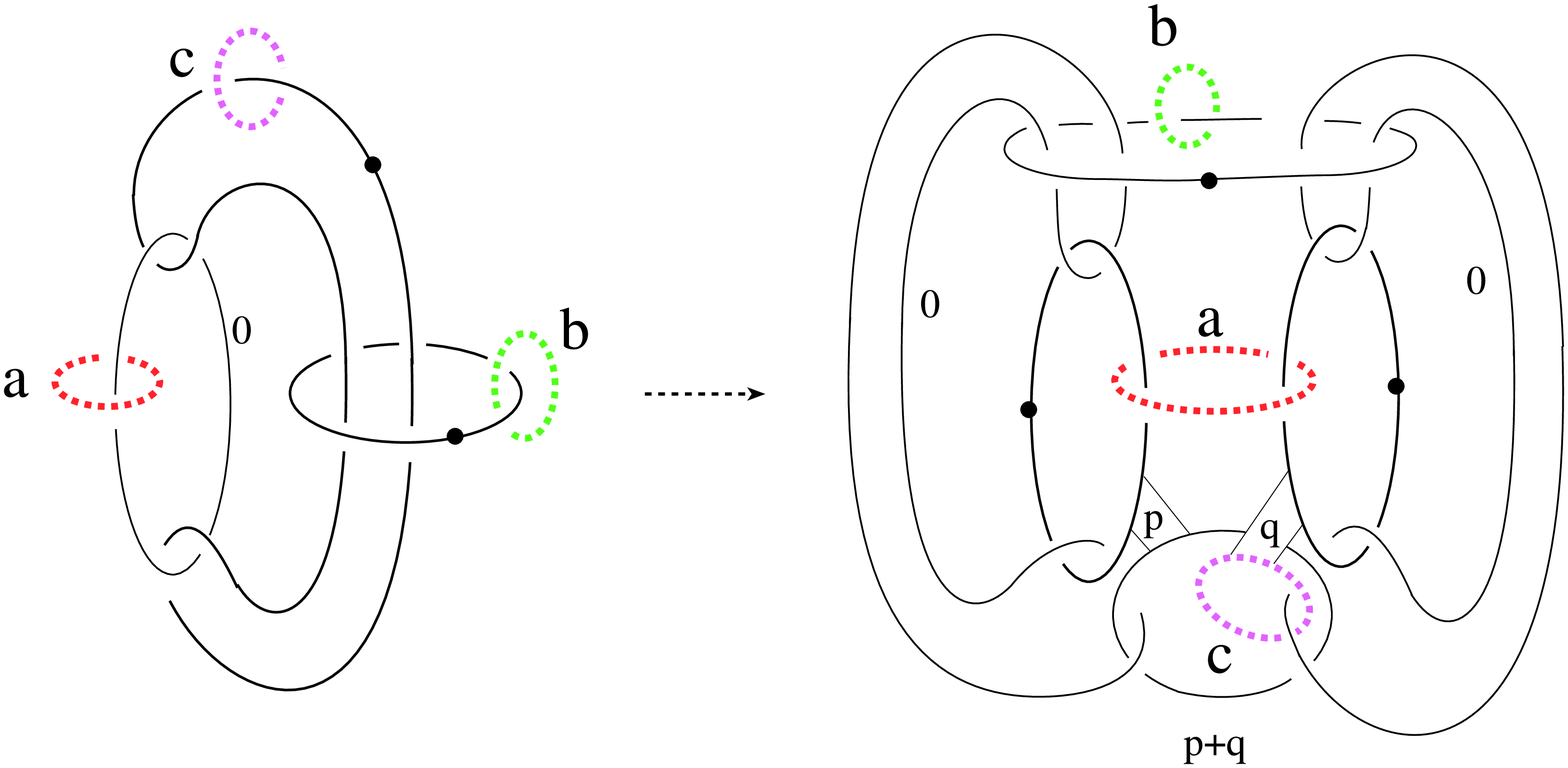}   
\caption{The algorithm $T^2\times B^2 \mapsto (T^2\times B^2)_{p,q} $} 
 \label{fig6}
\end{center}
\end{figure}

Given an imbedding $T^{2}\times B^2\subset X$, to describe the process  $X\mapsto X_{p,q}$ of doing p-log  and q-log transforms along two parallel tori, we need to visualize two distinct parallel copies of imbedded tori inside of $X$. To do this, in Figure~\ref{fig2} we start with the handlebody picture of $T^2$ ($2$-disk plus a pair of $1$-handles and a $2$-handle), then thicken this to $T^{2}\times [0,1]$, and then decompose it as two thickened tori identified along their faces $T^{2}\times [0,1/2]\cup T^{2}[1/2,1]$. Converting this description to Heegard diagram gives the first picture of Figure~\ref{fig3}, which is just  $T^{2}\times [0,1]$, where the two disjoint copies of $T^{2}$ inside are clearly visible, and then by thickening it we get $T^{2}\times B^2$ (the second picture of Figure~\ref{fig3}). Then by converting $1$-handle notation from pair of balls to ``circle with dot'' notation of \cite{a2}  we arrive to the pictures of $T^{2}\times B^2$ in Figure~\ref{fig4}, where the two disjoint copies of parallel $T^{2} \times B^{2}$ inside are visible. Now by using the recipe of Figure~\ref{fig1}, we perform p-log transform operation along one $T^{2}\times B^2$ , and the q-log operation along along the other copy of $T^{2}\times B^2$  we get Figure~\ref{fig5} which is $(T^{2}\times B^{2})_{p,q}$. Figure~\ref{fig6} summarizes the operation $X\mapsto X_{p,q}$. This figure indicates when you see a copy of $T^2\times B^2$ inside a handlebody $X^4$, how you perform the  (p,q)-log transform operation to $X$ (for this reason the images of the linking loops of $T^{2}\times B^2$ are indicated). In short, this algorithm is the combination of the process of pushing off two parallel copies of a given torus, followed by p-log transform operation of Figure~\ref{fig1}.

\vspace{.05in}

\begin{figure}[ht]  \begin{center}  
\includegraphics[width=.21\textwidth]{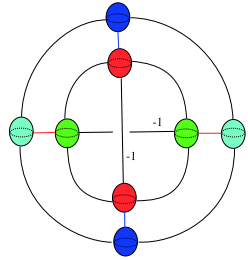}   
\caption{The cusp C}
 \label{fig7} 
\end{center}
\end{figure}

\begin{figure}[ht]  \begin{center}  
\includegraphics[width=.8\textwidth]{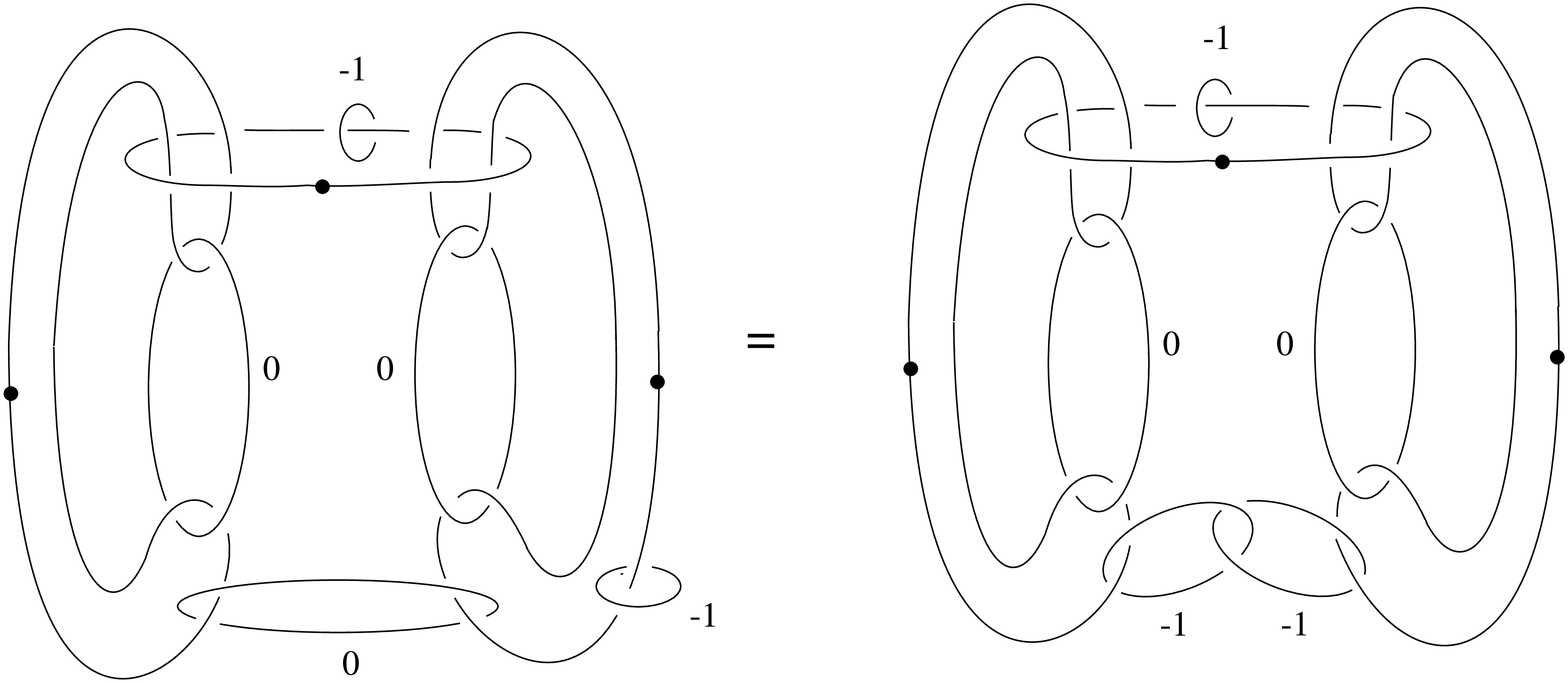}   
\caption{The cusp C in circle-with-dot notation} 
 \label{fig8}
\end{center}
\end{figure}

\begin{figure}[ht]  \begin{center}  
\includegraphics[width=.4\textwidth]{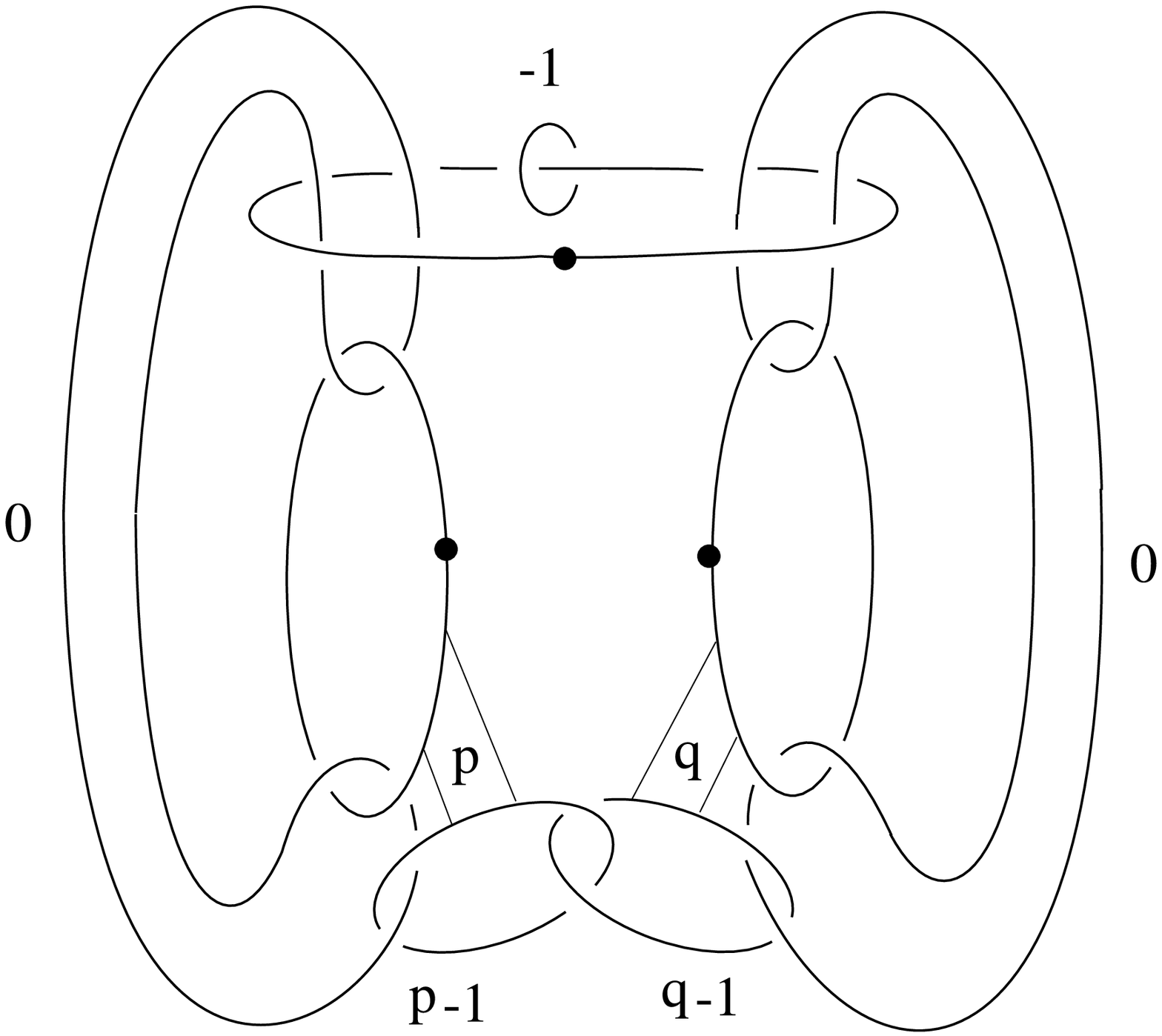}   
\caption{$C_{p,q}$} 
 \label{fig9}
\end{center}
\end{figure}

Now let us apply this operation to a cusp $C$, which is $T^2\times B^2$ plus two vanishing cycles (i.e. two $2$-handles) as shown in Figure~\ref{fig7}. By changing the 1-handle notation from pair of balls to circle with dot notation, and a handle slide, in Figure~\ref{fig8} we get two alternative pictures of C. Then finally, by applying the algorithm  above we get the handlebody picture of $C_{p,q}$ as shown in Figure~\ref{fig9}. 

\vspace{.05in}

\begin{figure}[ht]  \begin{center}  
\includegraphics[width=.74\textwidth]{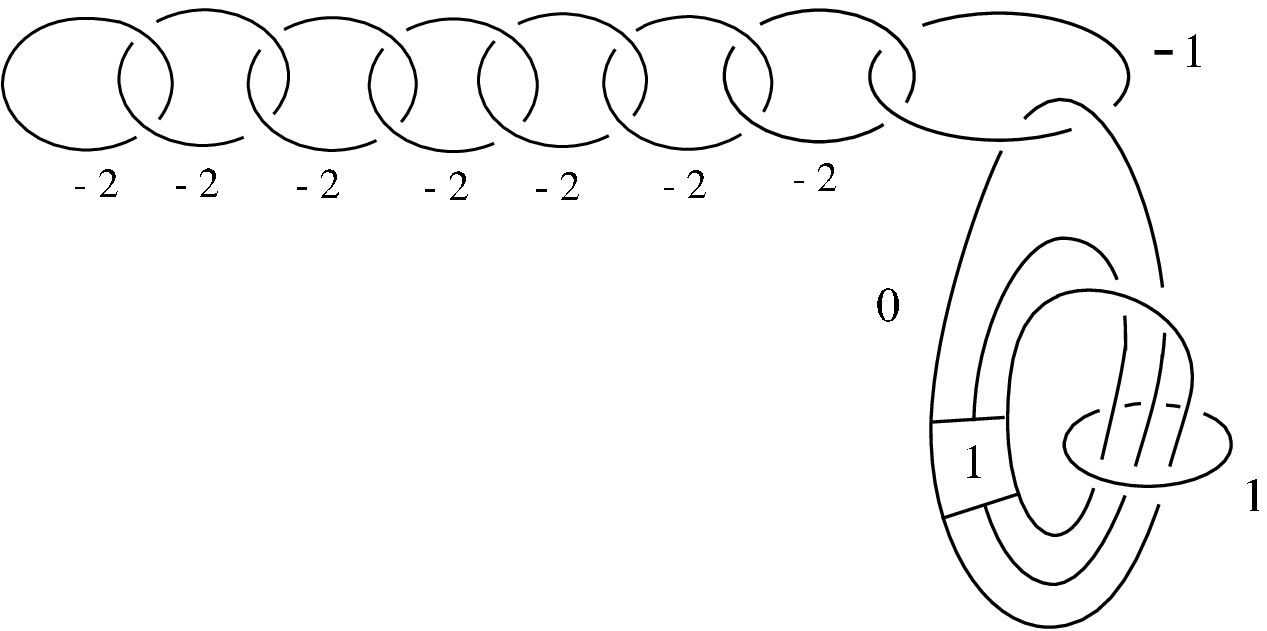}   
\caption{$E(1)$} 
 \label{fig10}
\end{center}
\end{figure}

\begin{figure}[ht]  \begin{center}  
\includegraphics[width=.5\textwidth]{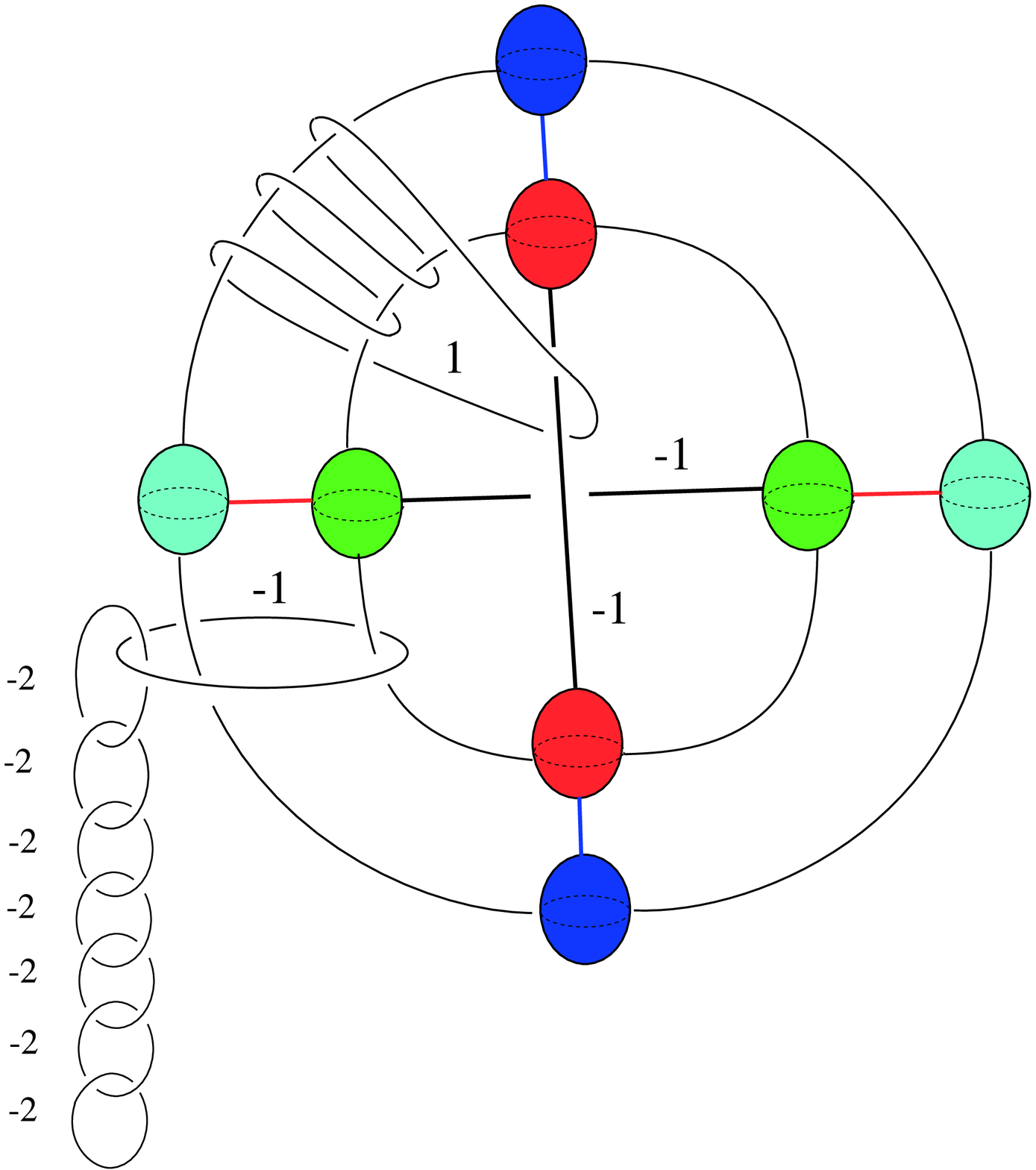}   
\caption{$E(1)$} 
 \label{fig11}
\end{center}
\end{figure}

We can easily extend  this process to $E(1)\mapsto E(1)_{p,q}$ as follows: First recall that the cusp $C$ sits inside of $E(1)$, as shown in the handlebody picture of $E(1)$ in Figure~\ref{fig10} (see \cite{a3}). By replacing the cusp $C$ inside with its equivalent  picture in Figure~\ref{fig7},  we get an another description of $E(1)$ in Figure~\ref{fig11} (the nice thing about this picture is that the two parallel copies of $T^2\times B^2$ inside of the cusp are clearly visible). Now by applying the process $C\mapsto C_{p,q}$ as before  (i.e. going through the steps Figure~\ref{fig7} $\mapsto$ Figure~\ref{fig9}), while carrying the rest of the handles of Figure~\ref{fig11}, we get a nice handlebody picture of $E(1)_{p,q}$ in Figure~\ref{fig12}.

\newpage

\section{Infinitely many exotic Stein manifolds}

  It is well known by a celebrated theorem of Eliashberg that, any smooth structure on the topological $B^4$ has to be diffeomorphic to the standard $B^4$, if it has a Stein structure. So in studying smooth structures on small  $4$-manifolds Stein structures play important role (one might think of them as force trying to keep the underlying smooth structure standard).
 In \cite{a4} an exotic pair of simple $4$-manifolds were introduced (each one is obtained from $B^4$ by attaching a $2$-handle). The main reason that they are exotic copies of each other is because one is a Stein manifold and the other is not.  After that, the question of whether exotic Stein manifold pairs can exist became important. Later on it turned out that they exist \cite{ay2}, and in fact a Stein manifold can have arbitrarily many exotic Stein copies of itself \cite{ay3}. ln particular, these examples can be made simply connected with  second Betti number $1$, and furthermore  their exoticness  can be explained by ``corks''. On the other hand surprisingly in  \cite{aems} infinitely many  exotic simply connected Stein manifolds with  large second Betti number were constructed. Here we find small (Betti number) such manifolds.

\vspace{.05in}

Now consider the handlebody of  $E(1)_{p,q}$ in Figure~\ref{fig12}, it contains the manifold $Q_{p,q}$ of Figure~\ref{fig13}  inside as a codimension zero submanifold. It follows from the construction that $\partial Q_{p,q}$ is the Poincare homology sphere $\Sigma=\Sigma(2,3,5)$. So we have a splitting  of $E(1)_{p,q}$  into two codimension zero pieces: $Q_{p,q}$ and its complement $N$, glued along $\Sigma$:
$$E(1)_{p,q}=Q_{p,q}\smile_{\Sigma} N$$ 
Clearly the piece $Q_{p,q}$ is responsible for the exoticity of $E(1)_{p,q}$.
From this decomposition we get yet another more interesting decomposition: 
$$E(1)_{p,q} \# \;5 \bar{\C\P}^2= M_{p,q}\smile_{\Sigma} N' $$
where $M_{p,q}$ is the manifold given in Figure~\ref{fig14}, and $N'$ is its complement.  This claim is evident from the imbedding $M_{p,q}\subset Q_{p,q} \# \;5 \bar{\C\P}^2$ (blowup $Q_{p,q}$ once on the $-1$ framed handle, and twice on each of  the $p-1$ and $q-1$ framed handles).  Also, by taking $p$ and $q$ relatively prime it is easy to see that  $M_{p,q}$ is simply connected, with second Betti number $3$. Moreover we claim  that $M_{p,q}$  is a Stein manifold. This can be seen by drawing the handlebody of Figure~\ref{fig14}  as a Legendrian handlebody of  Figure~\ref{fig15}, where each handle is $TB-1$ framed (a patient reader can check this by converting 1-handles of Figure~\ref{fig15} to ``circle with dot'' notation while keeping track of the framings). Hence we have:

\begin{Thm}  There are infinitely many simply connected, second Betti number $3$, Stein manifolds $M_{p,q}$, with a common boundary, which are exotic copies of each other relative to their boundaries. \end{Thm}

\section{Epilogue} 
\begin{figure}[ht]  \begin{center}  
\includegraphics[width=.6\textwidth]{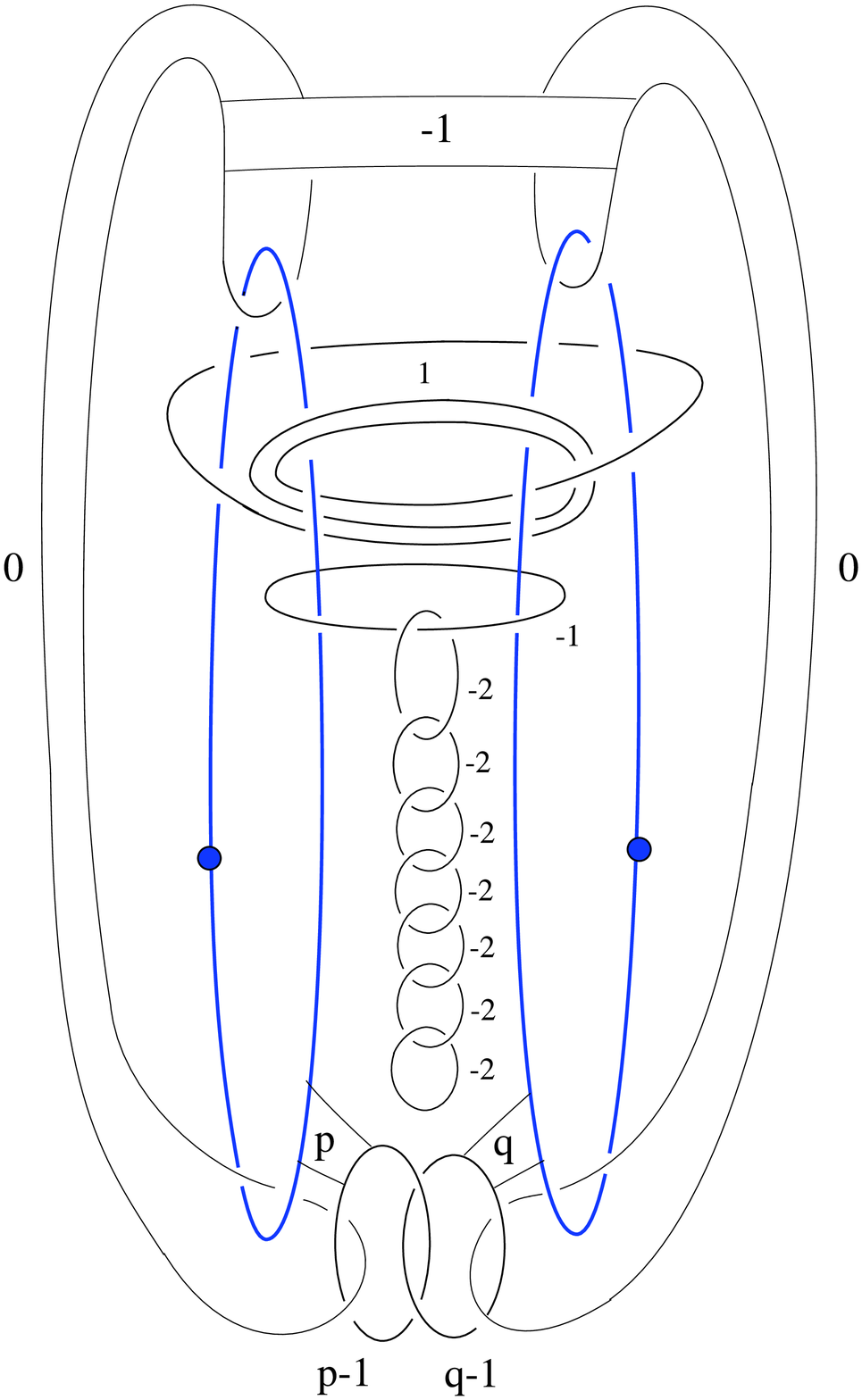}   
\caption{$E(1)_{p,q}$} 
 \label{fig12}
\end{center}
\end{figure}

\begin{figure}[ht]
\centering
\begin{minipage}{.5\textwidth}
  \centering
  \includegraphics[width=.7\linewidth]{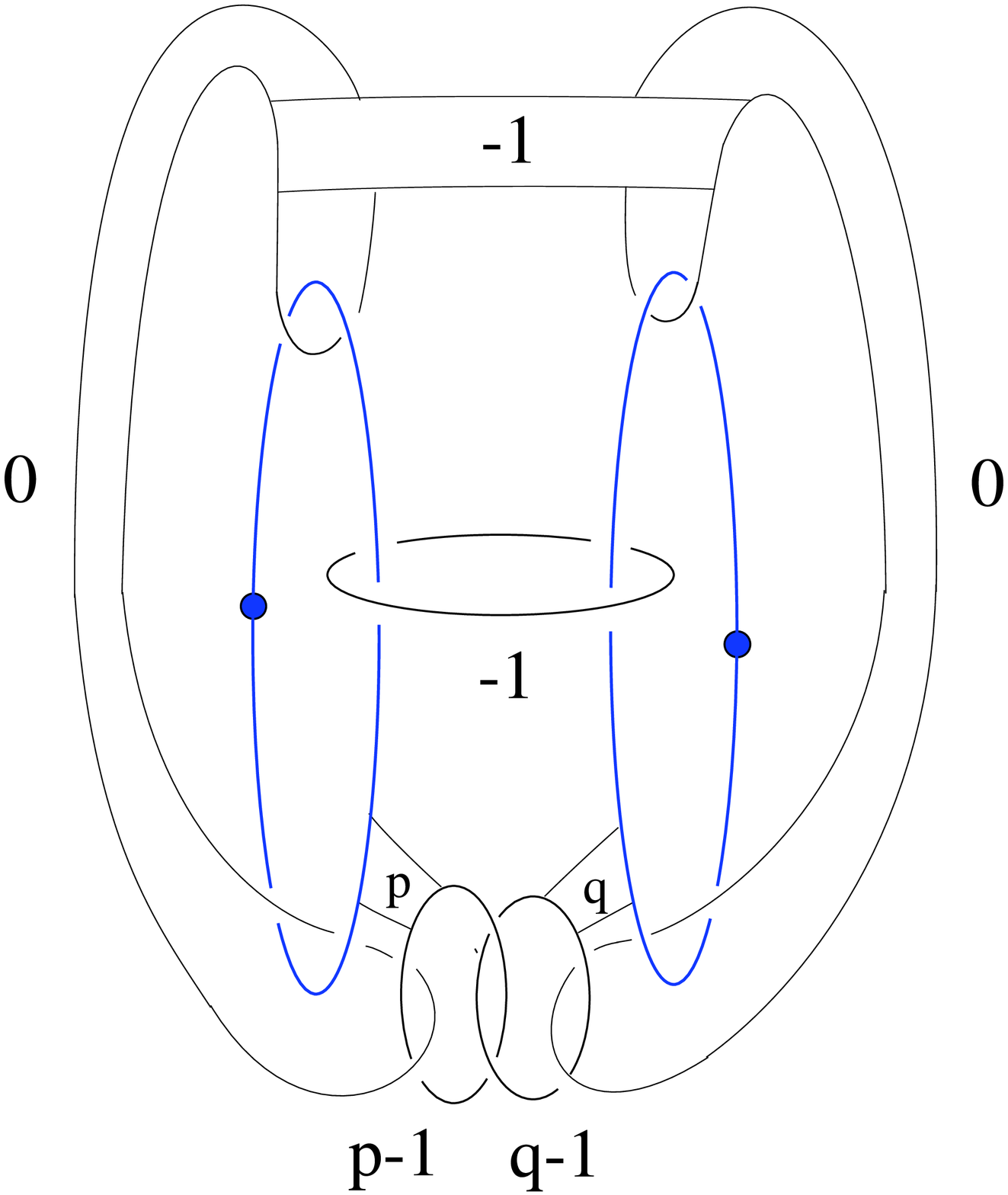}
  \caption{$Q_{p,q}$}
  \label{fig13}
\end{minipage}%
 \begin{minipage}{.5\textwidth}
  \includegraphics[width=.7\linewidth]{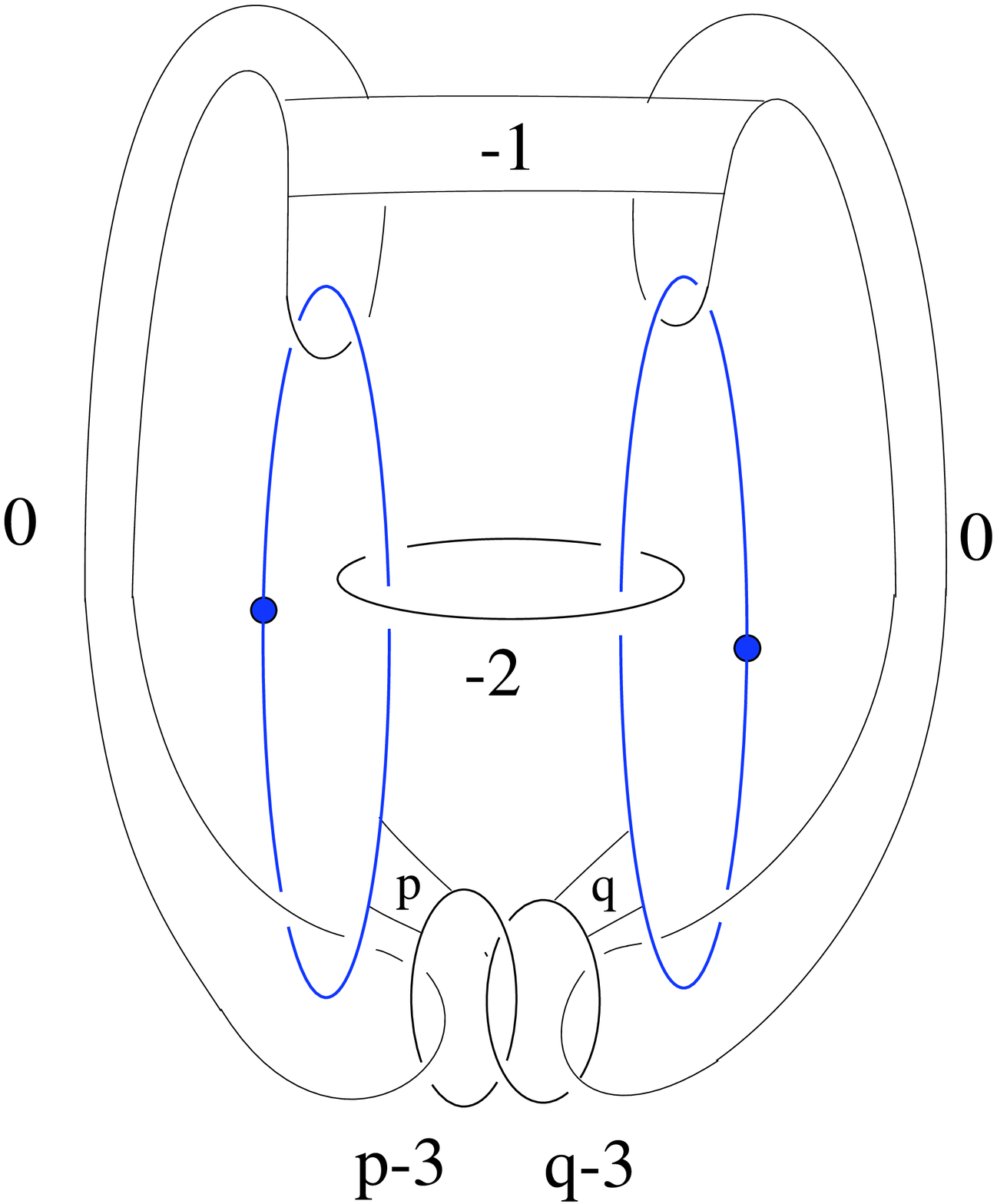}
  \caption{$M_{p,q}$}
  \label{fig14}
\end{minipage}
\end{figure}

\begin{figure}[ht]  \begin{center}  
\includegraphics[width=.7\textwidth]{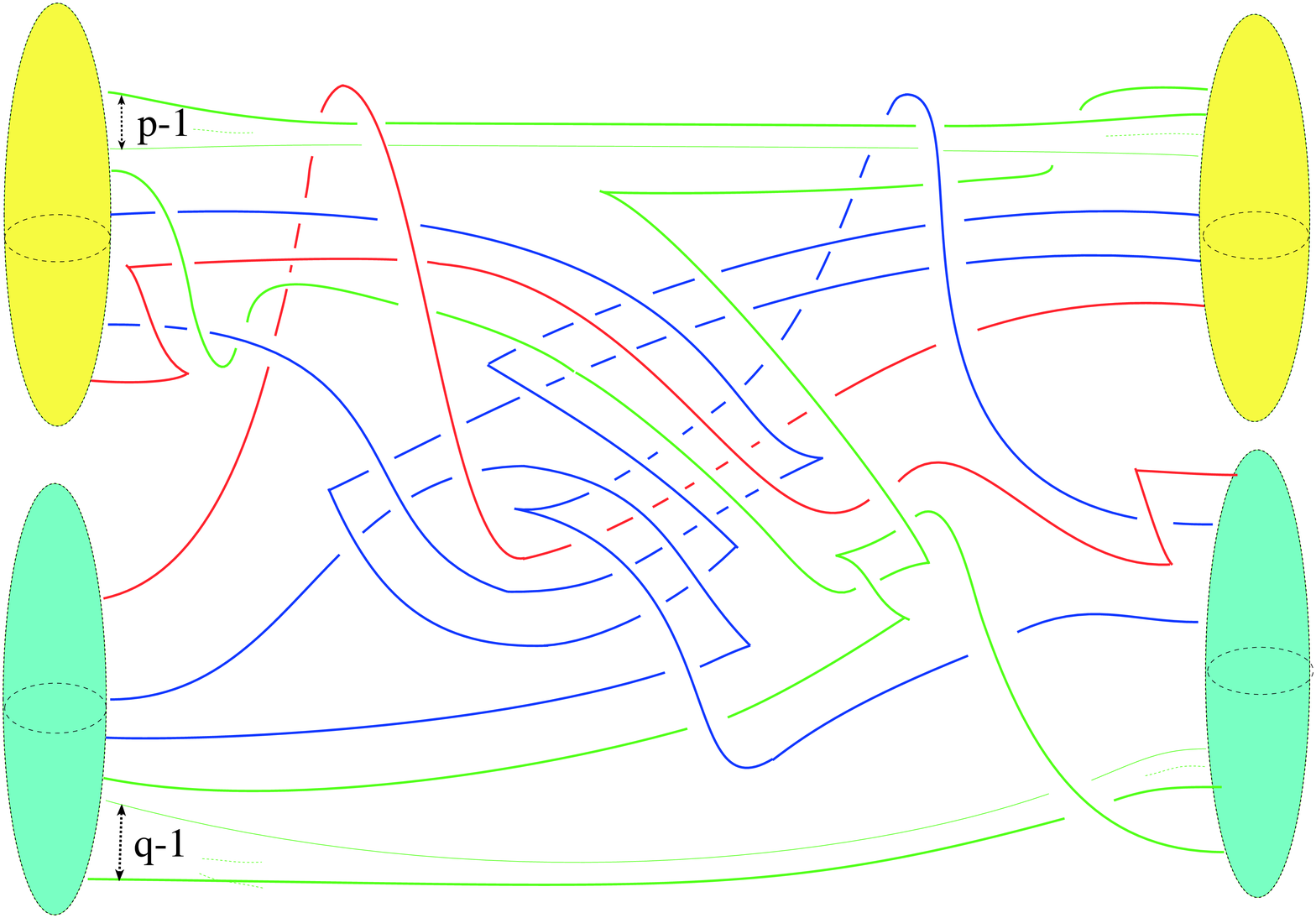}   
\caption{$M_{p,q}$} 
 \label{fig15}
\end{center}
\end{figure}

After first posting of this paper, the author and K. Yasui have found infinitely many exotic Stein fillings $X^{4}_{p}$ of a fixed contact $3$-manifold $Y^3$ \cite{ay4}. The manifolds $X_{p}$ are all simply connected and has Betti number $2$, and they are exotic copies of each other (without the ``relative to boundary'' condition of Theorem 1). They seem to be related to the Stein manifolds $M_{p,q}$ appearing inside of the blown-up copies of $E(1)_{p,q}$. For completeness, here we make this relation  concrete. In  Figure~\ref{fig12} by sliding $q-1$ framed $2$-handle over the $-1$ framed $2$-handle $q$ times, we make it disjoint from the $1$-handle which  it links originally. This gives the manifold $Z'_{p,q}\subset E(1)_{p,q}$ of  Figure~\ref{fig20}. By decreasing the framings of the $2$-handles of $Z_{p,q}'$ by one,  we obtain an imbedding of a Stein manifold  $Z_{p,q} \subset E(1)_{p,q}\# 2\bar{\C\P}^2 $ of Figure~\ref{fig21} (this step is similar to Figure~\ref{fig13} $\mapsto $  Figure~\ref{fig14}). We drew the $1$-handles of this manifold in pair of balls notation, to exhibit its Stein structure. By using \cite{b}, it is easy to check that, the diffeomorphism $f: \partial (Z_{p,q})\to \partial (Z_{p+2,q})$ induced from log-p transform operation (applied twice)  extends to a homeomorphism inside. 

\vspace{.05in}

When $q=2$ and $p$ odd, the Stein manifolds  $Z_{p}:= Z_{p,2}$ are simply connected and has second Betti numbers equal to $2$,  and $H_{2}(Z_{p})$ is generated by the classes $C$ and $D = 2A-pB+(p^{2}+1)C$ (where $A, B, C$ are the $2$-handles of Figure~\ref{fig21}) with the intersection form: 

$$\left( \begin{array}{cc}
0 & 2  \\
2 & -4 \end{array} \right) $$

It follows that any self intersection $-4$ homology class has to be $\pm D$. Let $g(D)$ be the genus of any surface representing $D$. As in \cite{ay4}, by applying  the version of the adjunction inequality in \cite{am} we get:
$$ 2g(D)-2 \geq D.D + pq=-4 + pq$$

\noindent (from Figure~\ref{fig21} the rotational numbers of A, B, C can be readily computed to be  0, q, 0,  respectively). Hence the minimal genus of the homology class $D$ approaches infinity, $g(D) \to \infty $ as $p\to \infty$. Therefore  the sequence of $4$-manifolds $Z_{p}$ as $p\to \infty$ contains infinitely many distinct exotic copies of each other.

\newpage

\section{$(S^{2}\times S^{2})_{p,q} =S^{2}\times S^{2}$ }

It is known that  $S^2\times S^2$ is a union of two cusps $C\cup_{id} -C$ glued along their boundaries, as shown in the first picture of Figure~\ref{fig17}.  So it is natural to
ask whether we get exotic copies of  $S^2\times S^2$  by performing (p,q)-log transforms to $S^2\times S^2$ along the parallel tori in these cusps. Here we answer this question negatively, by proving a result similar to the one proved in \cite{a5} for the case of knot surgeries. 
First of all since the tori inside of $-C$  isotope to the tori inside of $C$, it suffices to prove $(S^2\times S^2)_{p,q}:=C_{p,q}\cup -C$ is standard copy of $S^2\times S^2$. 
\begin{Thm}  $(S^2\times S^2)_{p,q}$ is diffeomorphic to $S^2\times S^2$. 
\end{Thm}
\proof We start with the handlebody $S^2\times S^2 =C\cup_{id} -C$ which is the first picture of Figure~\ref{fig17}.  By using our algorithm we first perform the operation $S^2\times S^2 \mapsto (S^2\times S^2)_{p,q}$ we get the second picture, then by canceling the  middle $0$-framed $2$-handle by one of the $1$-handles we get the third picture of Figure~\ref{fig17}. Then after a handle slide (as indicated by the dotted arrow) we get the last picture of Figure~\ref{fig17}. Then by another handle slide (indicated by the dotted arrow) we arrive to the first picture of Figure~\ref{fig18}. By further handle slides, and  by introducing a canceling pair of $2$- and $3$ handles, then a handle slide followed by taking away canceling $2$- and $3$- handle pairs we obtain the last picture of Figure~\ref{fig18} (these steps are similar to the one used in \cite{a6}). Then by an isotopy and handle slides we arrive to the first picture of Figure~\ref{fig19} ($m=p+q$). Evidently when $m$ and $q$ are relatively prime, this is $S^2\times S^2$ (to see this, keep sliding the two $2$-handles, which link the $1$- handle, over each other, and use the fact that the excess (non-algebraic) geometric linkings can be removed by the help of the small  $0$-framed dual handles). Then after canceling $1/2$ and $2/3$ handle pairs we arrive to the middle then to the last picture of Figure~\ref{fig19}, which is $S^{2}\times S^{2}$. \qed

\newpage

 \begin{figure}[ht]  \begin{center}  
\includegraphics[width=.6\textwidth]{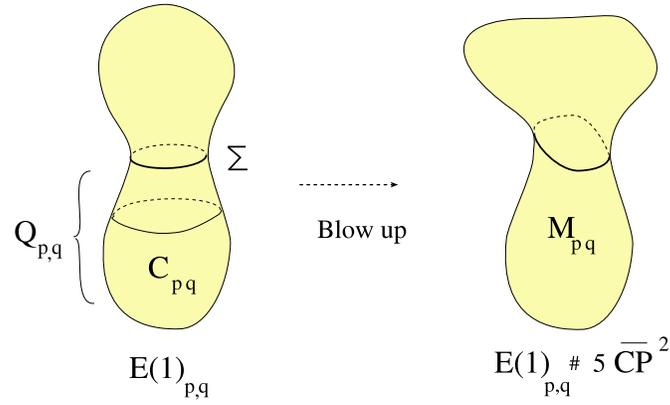}   
\caption{Decompositions of $E(1)_{p,q}$ and $E(1)_{p,q}  \# 5 \bar{\C\P^{2}}$ } 
 \label{fig16}
\end{center}
\end{figure}

\begin{figure}
\centering
\begin{minipage}{.4\textwidth}
  \centering
  \includegraphics[width=.8\linewidth]{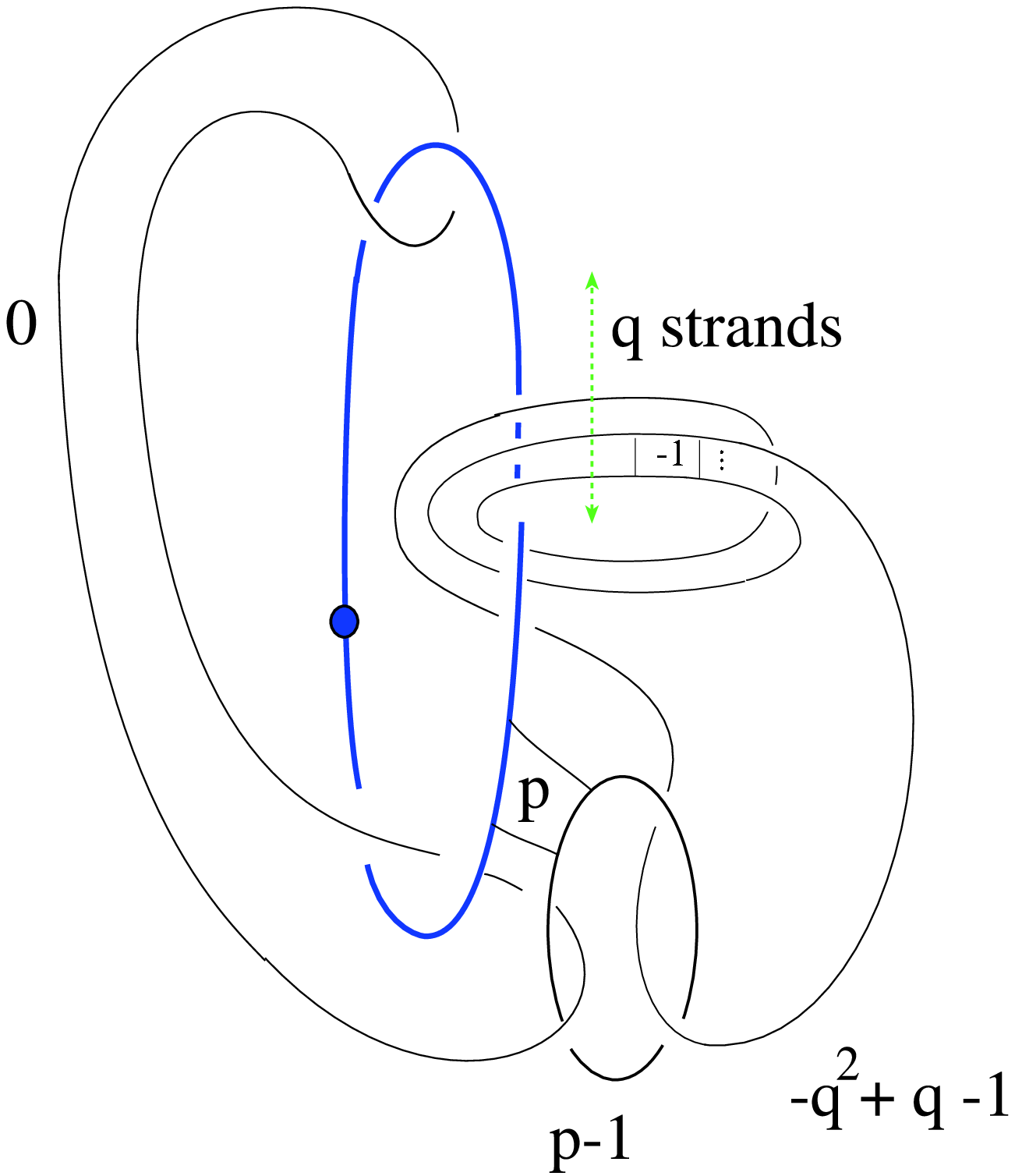}
  \captionof{figure}{$Z'_{p,q}$}
  \label{fig20}
\end{minipage}%
\begin{minipage}{.5\textwidth}
  \centering
  \includegraphics[width=.9\linewidth]{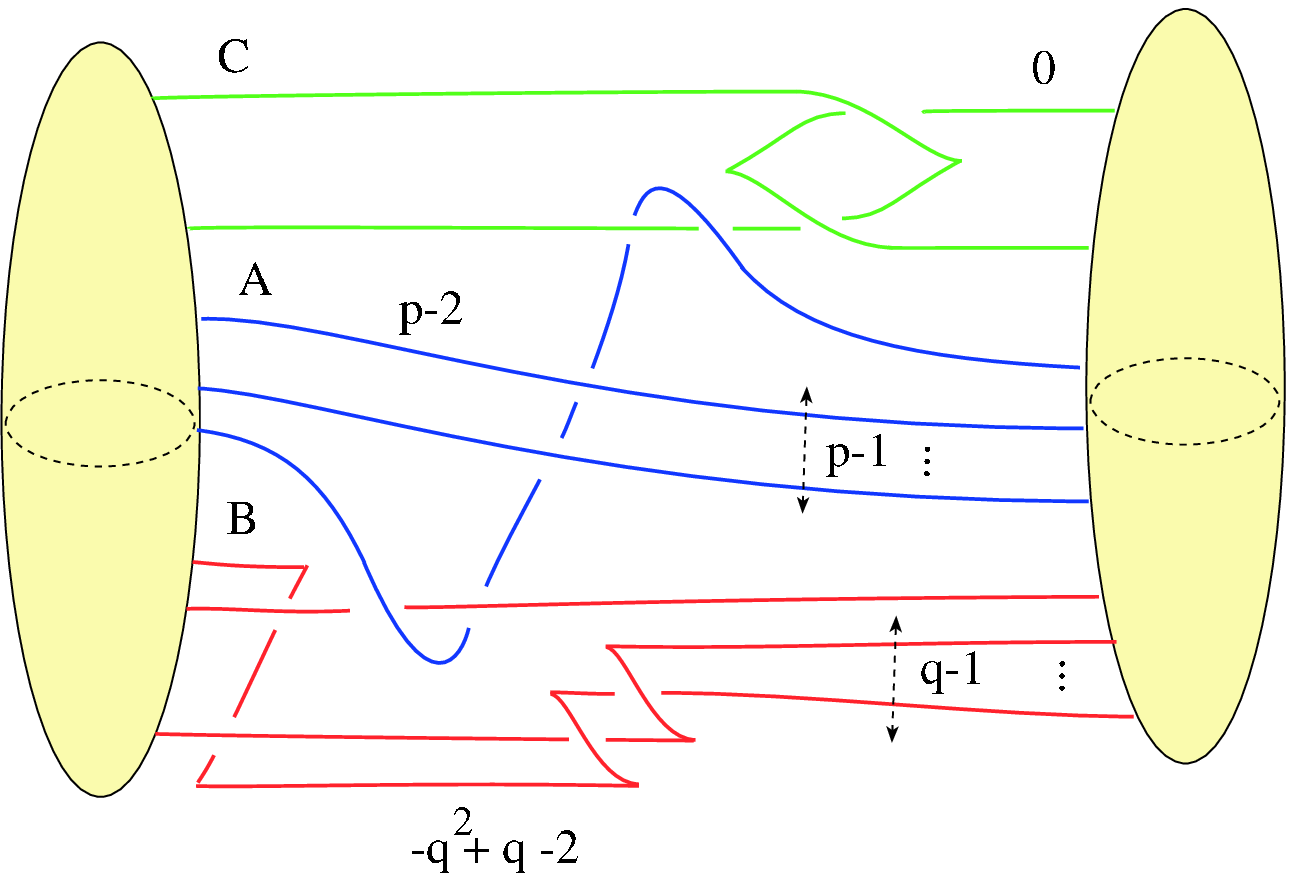}
  \captionof{figure}{$Z_{p,q}$}
  \label{fig21}
\end{minipage}
\end{figure}

  \begin{figure}[ht]  \begin{center}  
\includegraphics[width=.8\textwidth]{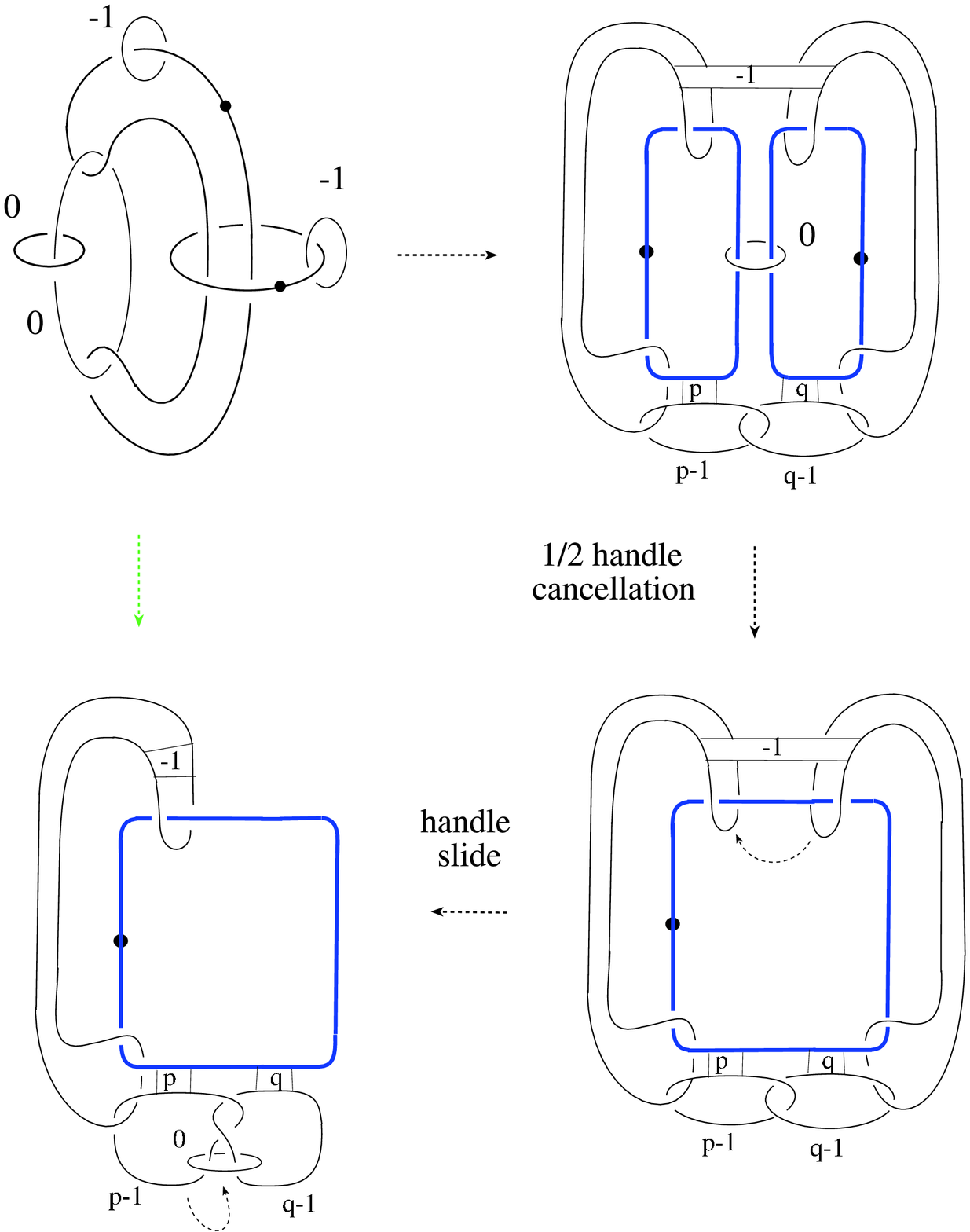}   
\caption{$S^2\times S^2=C\cup -C \mapsto (S^2\times S^2)_{p,q}$} 
 \label{fig17}
\end{center}
\end{figure} 

 \begin{figure}[ht]  \begin{center}  
\includegraphics[width=.8\textwidth]{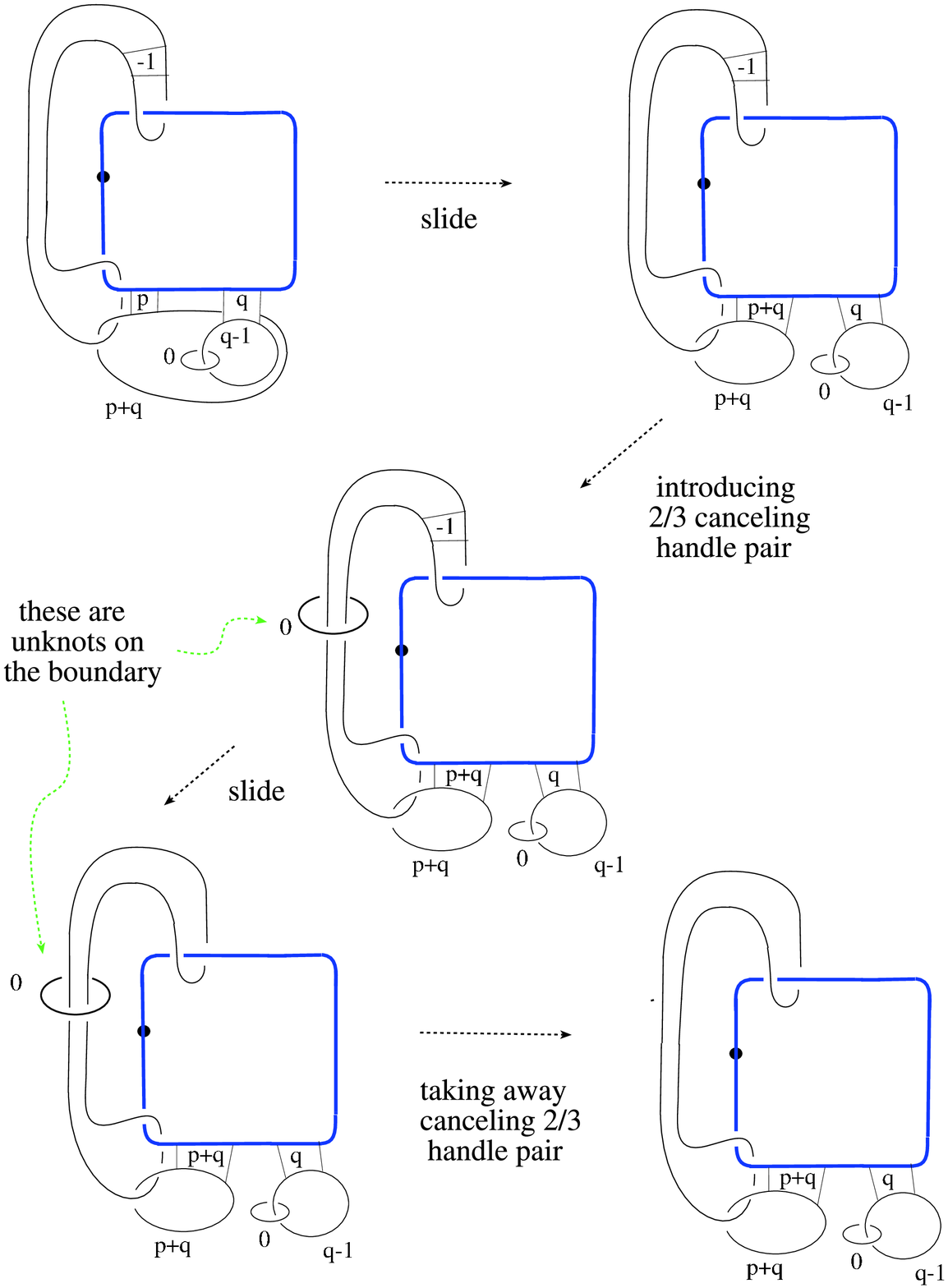}   
\caption{$(S^2\times S^2)_{p,q} \mapsto S^2\times S^2$} 
 \label{fig18}
\end{center}
\end{figure} 

\clearpage

  \begin{figure}[ht]  \begin{center}  
\includegraphics[width=.7\textwidth]{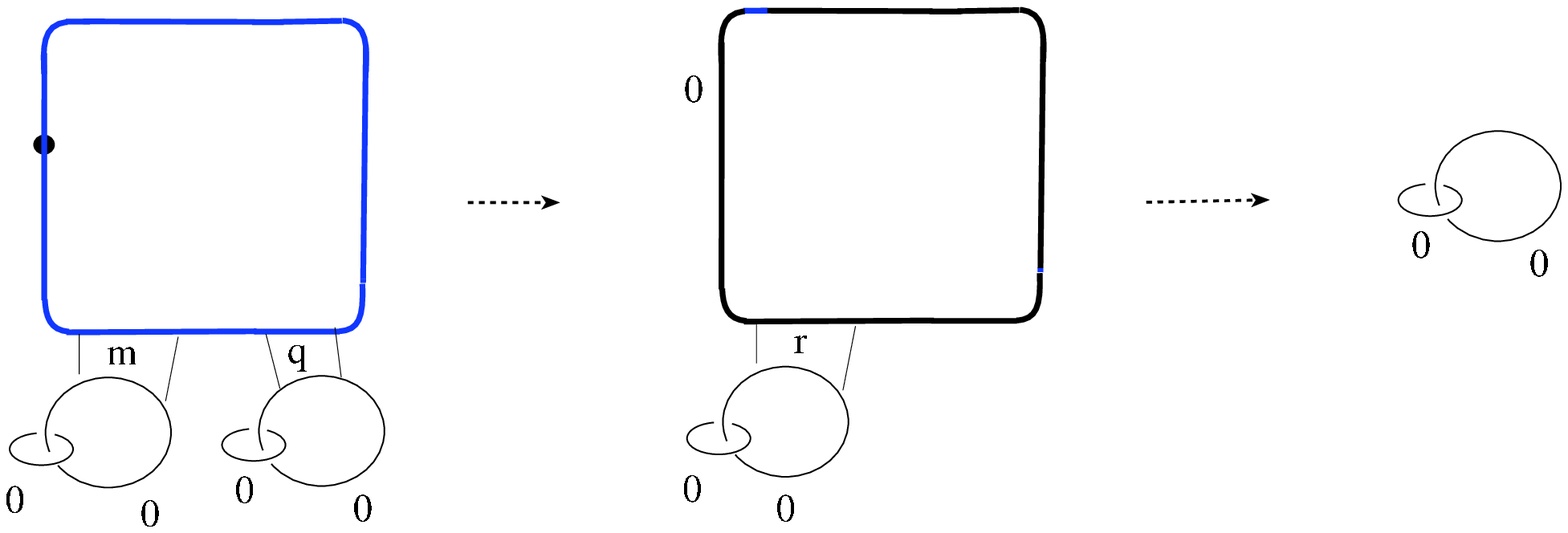}   
\caption{$(S^2\times S^2)_{p,q} \mapsto S^2\times S^2$ continued} 
 \label{fig19}
\end{center}
\end{figure}

\end{document}